\newcommand{\CLASSINPUTtoptextmargin}
{0.78in}

\newcommand{\CLASSINPUTbottomtextmargin}
{1.02in}

\documentclass[conference]{IEEEtran}

\makeatletter
\def\ps@headings{%
\def\@oddhead{\mbox{}\scriptsize\rightmark \hfil \thepage}%
\def\@evenhead{\scriptsize\thepage \hfil \leftmark\mbox{}}%
\def\@oddfoot{}%
\def\@evenfoot{}}
\makeatother 

\usepackage{cite}      
\usepackage{subfigure} 
\usepackage{url}       
\usepackage{algorithm}
\usepackage{array}
\usepackage{amsmath,amsfonts,amssymb,graphicx}
\usepackage{bm}
\usepackage{mathrsfs}
\usepackage{dsfont}
\usepackage{multicol}
\usepackage[end]{algpseudocode}
\usepackage{multirow}

\newtheorem{theorem}{\textbf{Theorem}}
\newtheorem{lemma}{\textbf{Lemma}}

\newcommand {\xA} {\mathcal{A}}

\newcommand{\bX} {\hat{\overline{\theta}}}
\newcommand{\obX} {\hat{\overline{\theta^*}}}

\newcommand{\hX} {\hat{\theta}}

\newcommand{\tT} {\widetilde{T}_{i,j}}
\newcommand{\tI} {\widetilde{I}_{i,j}}

\newcommand{\commnt}[1] {$//$ \textsc{#1} }
\newcounter{MYtempeqncnt}

\begin{document}
\title{On the Combinatorial Multi-Armed Bandit Problem with Markovian Rewards}

\author{\IEEEauthorblockN{Yi Gai$^*$, Bhaskar Krishnamachari$^*$ and Mingyan Liu$^\ddag$}
\IEEEauthorblockA{$^*$Ming Hsieh Department of Electrical Engineering\\
University of Southern California, Los Angeles, CA 90089, USA\\
$^\ddag$Department of Electrical Engineering and Computer Science\\
University of Michigan, Ann Arbor, MI 48109, USA\\
 Email: $\{$ygai,bkrishna$\}$@usc.edu; mingyan@eecs.umich.edu}}

\maketitle

\begin{abstract}

We consider a combinatorial generalization of the classical
multi-armed bandit problem that is defined as follows. There is a
given bipartite graph of $M$ users and $N \geq M$ resources. For
each user-resource pair $(i,j)$, there is an associated state that
evolves as an aperiodic irreducible finite-state Markov chain with
unknown parameters, with transitions occurring each time the
particular user $i$ is allocated resource $j$. The user $i$ receives
a reward that depends on the corresponding state each time it is
allocated the resource $j$. The system objective is to learn the
best matching of users to resources so that the long-term sum of the
rewards received by all users is maximized. This corresponds to
minimizing regret, defined here as the gap between the expected
total reward that can be obtained by the best-possible static
matching and the expected total reward that can be achieved by a
given algorithm. We present a polynomial-storage and
polynomial-complexity-per-step matching-learning algorithm for this
problem. We show that this algorithm can achieve a regret that is
uniformly arbitrarily close to logarithmic in time and polynomial in
the number of users and resources. This formulation is broadly
applicable to scheduling and switching problems in networks and
significantly extends prior results in the area.
\end{abstract}

\section{Introduction}\label{sec:intro}

Multi-armed bandit problems provide a fundamental approach to
learning under stochastic rewards, and find rich applications in a
wide range of networking contexts, from Internet
advertising~\cite{Pandey} to medium access in cognitive radio
networks~\cite{Gai:2010, Liu:Zhao, Anandkumar:2010}. In the
simplest, classic non-Bayesian version of the problem, studied by
Lai and Robbins~\cite{Lai:Robbins}, there are K independent arms,
each generating stochastic rewards that are i.i.d. over time. The
player is unaware of the parameters for each arm, and must use some
policy to play the arms in such a way as to maximize the cumulative
expected reward over the long term. The policy's performance is
measured in terms of its ``regret", defined as the gap between the
the expected reward that could be obtained by an omniscient user
that knows the parameters for the stochastic rewards generated by
each arm and the expected cumulative reward of that policy. It is of
interest to characterize the growth of regret with respect to time
as well as with respect to the number of arms/players. Intuitively,
if the regret grows sublinearly over time, the time-averaged regret
tends to zero.

There is inherently a tradeoff between exploration and exploitation
in the learning process in a multi-armed bandit problem: on the one
hand all arms need to be sampled periodically by the policy used, to
ensure that the "true" best arm is found; on the other hand, the
policy should play the arm that is considered to be the best often
enough to accumulate rewards at a good pace.

In this paper, we formulate a novel combinatorial generalization of
the multi-armed bandit problem that allows for Markovian rewards and
propose an efficient policy for it. In particular, there is a given
bipartite graph of $M$ users and $N \geq M$ resources. For each
user-resource pair $(i,j)$, there is an associated state that
evolves as an aperiodic irreducible finite-state Markov chain with
unknown parameters, with transitions occurring each time the
particular user $i$ is allocated resource $j$. The user $i$ receives
a reward that depends on the corresponding state each time it is
allocated the resource $j$. A key difference from the classic
multi-armed bandit is that each user can potentially see a different
reward process for the same resource. If we therefore view each
possible matching of users to resources as an arm, then we have a
super-exponential number of arms with dependent rewards. Thus, this
new formulation is significantly more challenging than the
traditional multi-armed bandit problems.

Because our formulation allows for user-resource matching, it could
be potentially applied to a diverse range of networking settings
such as switching in routers (where inputs need to be matched to
outputs) or frequency scheduling in wireless networks (where nodes
need to be allocated to channels) or for server assignment problems
(for allocating computational resources for various processes),
etc., with the objective of learning as quickly as possible so as to
maximize the usage of the best options. For instance, our
formulation is general enough to be applied to the channel
allocation problem in cognitive radio networks considered
in~\cite{Gai:2010} if the rewards for each user-channel pair come
from a discrete set and are i.i.d. over time (which is a special
case of Markovian rewards).

Our main contribution in this work is the design of a novel policy
for this problem that we refer to Matching Learning for Markovian
Rewards (MLMR). Since we treat each possible matching of users to
resources as an arm, the number of arms in our formulation grows
super-exponentially. However, MLMR uses only polynomial storage, and
requires only polynomial computation at each step. We analyze the
regret for this policy with respect to the best possible static
matching, and show that it is uniformly logarithmic over time under
some restrictions on the underlying Markov process. We also show
that when these restrictions are removed, the regret can still be
made arbitrarily close to logarithmic with respect to time. In
either case, the regret is polynomial in the number of users and
resources.

The rest of the paper is organized as follows. In
section~\ref{sec:priorwork} we present our work in the context of
prior results on multi-armed bandits. In
section~\ref{sec:formulation} we present the problem formulation. In
section~\ref{sec:algorithm} we present a polynomial-storage
polynomial-time-per-step learning policy, which we refer to as MLMR.
We analyze the regret for this policy in section~\ref{sec:regret}
and show that it yields a bound on the regret that is uniformly
logarithmic over time and polynomial in the number of users and
resources under certain conditions on the Markov chains describing
the state evolution for the arms. We then show that the regret can
still be arbitrarily close to logarithmic with respect to time when
no knowledge is available. We present some examples and simulations
in section~\ref{sec:simulation}, and conclude with some comments and
ideas for future work in section~\ref{sec:conclusion}.

\section{Prior Work}\label{sec:priorwork}

The problem we consider in this paper is different from prior work
for two key reasons. We treat rewards that are dependent across a
super-exponential number of arms whose states evolve in a non-i.i.d.
Markovian fashion over time. We summarize below prior work, which
has treated a) independent and temporally i.i.d. rewards, or b)
independent and Markovian state-based rewards, or c) non-independent
arms with temporally i.i.d.

\subsection{Independent arms with temporally i.i.d. rewards}
The work by Lai and Robbins~\cite{Lai:Robbins} assumes K independent
arms, each generating rewards that are i.i.d. over time from a given
family of distributions with an unknown real-valued parameter. For
this problem, they present a policy that provides an expected regret
that is $O(K \log n)$, i.e. linear in the number of arms and
asymptotically logarithmic in n. Anantharam \emph{et al.} extend
this work to the case when $M$ simultaneous plays are
allowed~\cite{Anantharam}. The work by Agrawal~\cite{Agrawal:1995}
presents easier to compute policies based on the sample mean that
also has asymptotically logarithmic regret. The paper by Auer
\emph{et al.}~\cite{Auer:2002} that considers arms with non-negative
rewards that are i.i.d. over time with an arbitrary un-parameterized
distribution that has the only restriction that it have a finite
support. Further, they provide a simple policy (referred to as
UCB1), which achieves logarithmic regret uniformly over time, rather
than only asymptotically. Our work utilizes a general
Chernoff-Hoeffding-bound-based approach to regret analysis pioneered
by Auer \emph{et al.}.

Some recent work has shown the design of distributed multiuser
policies providing asymptotically logarithmic regret, for the
context of cognitive radio networks~\cite{Liu:Zhao,
Anandkumar:2010}.

\subsection{Independent arms with Markovian rewards}

There has been relatively less work on multi-armed bandits with
Markovian rewards. Anantharam \emph{et al.}~\cite{Anantharam:1987}
wrote one of the earliest papers with such a setting. They proposed
a policy to pick $m$ out of the $N$ arms each time slot and prove
the lower bound and the upper bound on regret. However, the rewards
in their work are assumed to be generated by rested Markov chains
with transition probability matrices defined by a single parameter
$\theta$ with identical state spaces. Also, the result for the upper
bound is achieved only asymptotically.

For the case of single users and independent arms, a recent work by
Tekin and Liu~\cite{Tekin:2010} has extended the results in
\cite{Anantharam:1987} to the case with no requirement for a single
parameter and identical state spaces across arms. They propose to
use UCB1 from \cite{Auer:2002} for the multi-armed bandit problem
with rested Markovian rewards and prove a logarithmic upper bound on
the regret under some conditions on the Markov chain. We use
elements of the proof from~\cite{Tekin:2010} in this work, which is
however quite different in its combinatorial matching formulation
(which allows for dependent arms). The work on restless Markovian
rewards with single users and independent arms could be found in
\cite{Tekin:restless, Qing:restless, Dai:restless}.

\subsection{Dependent arms with temporally i.i.d. rewards}

The paper by Pandey \emph{et al.}~\cite{Pandey} divides arms into
clusters of dependent arms, each providing binary rewards, but they
do not present any theoretical analysis on the expected regret.
In~\cite{Rusmevichientong}, the reward from each arm is modeled as
the sum of a linear combination of a set of static random numbers
and a zero-mean random variable that is i.i.d. over time and
independent across arms. This is quite different from our model of
rewards.

Our work in this paper is closest to and builds on the recent work
which introduced combinatorial multi-armed bandits~\cite{Gai:2010}.
The formulation in~\cite{Gai:2010} has the restriction that the
reward process must be i.i.d. over time. A polynomial storage
matching learning algorithm is presented in~\cite{Gai:2010} that
yields regret that is polynomial in users and resources and
uniformly logarithmic in time for the case of i.i.d. rewards.
Although i.i.d. rewards are a special case of Markovian state-based
rewards, one reason this work is not a strict generalization of
~\cite{Gai:2010} is our assumption that the number of possible
states, and hence the support of the reward distribution on each
arm, is finite (whereas \cite{Gai:2010} allows for continuous reward
distributions with bounded support).

\section{Problem Formulation}\label{sec:formulation}

We consider a bipartite graph with $M$ users and $N \geq M$
resources predefined by some application. Time is slotted and is
indexed by $n$. At each decision period (also referred to
interchangeably as time slot), each of the $M$ users is assigned a
resource with some policy.

For each user-resource pair $(i,j)$, there is an associated state
that evolves as an aperiodic irreducible finite-state Markov chain
with unknown parameters. When user $i$ is assigned resource $j$,
assuming there are no other conflicting users assigned this
resource, $i$ is able to receive a reward that depends on the
corresponding state each time it is allocated the resource $j$.  We
denote the state space as $S_{i,j} = \{z_1, z_2, \ldots,
z_{|S_{i,j}|}\}$. The state of the Markov chain for each
user-resource pair $(i,j)$ evolves only when resource $j$ is
allocated to user $i$. We assume the Markov chains for different
user-resource pairs are mutually independent. The reward got by user
$i$ while allocated resource $j$ on state $z \in S_{i,j}$ is denoted
as $\theta^{i,j}_z$, which is also unknown to the users. We denote
$\mathbf{P}_{i,j} = \{p_{i,j}(z_a,z_b)\}_{z_a, z_b \in S_{i,j}}$ as
the transition probability matrix for the Markov chain $(i,j)$.
Denote $\pi^{i,j}_z$ as the steady state distribution for state $z$.
The mean reward got by user $i$ on resource $j$ is denoted as
$\mu_{i,j}$. Then we have $\mu_{i,j} = \sum\limits_{z \in S_{i,j}}
\theta^{i,j}_z \pi^{i,j}_z$. The set of all mean rewards is denoted
as $\bm{\mu} = \{\mu_{i,j}\}$.

We denote $Y_{i,j}(n)$ as the actual reward obtained by a user $i$
if it is assigned resource $j$ at time $n$. We assume that
$Y_{i,j}(n) = \theta_{i,j}^{z(n)}$, if user $i$ is the only occupant
of resource $j$ at time $n$ where $z(n)$ is the state of Markov
chain associated with $(i,j)$ at time $n$. Else, if multiple users
are allocated resource $j$, then we assume that, due to
interference, at most one of the conflicting users $j'$ gets reward
$Y_{i,j'}(n) = \theta_{i,j'}^{z'(n)}$ where $z'(n)$ is the state of
Markov chain associated with $(i,j')$ at time $n$,  while the other
users on the resources $j \neq j'$ get zero reward, i.e.,
$Y_{i,j}(n) = 0$. This interference model covers scenarios in many
networking settings.

A deterministic policy $\alpha(n)$ at each time is defined as a map
from the observation history $\{\mathbf{O}_t\}_{t=1}^{n-1}$ to a
vector of resources $\mathbf{o}(n)$ to be selected at period $n$,
where $\mathbf{O}_t$ is the observation at time $t$; the $i$-th
element in $\mathbf{o}(n)$, $o_i(n)$, represents the resource
allocation for user $i$. Then the observation history
$\{\mathbf{O}_t\}_{t=1}^{n-1}$ in turn can be expressed as
$\{o_i(t), Y_{i, o_i(t)}(t)\}_{1 \leq i \leq M, 1 \leq t < n}$.

Due to the fact that allocating more than one user to a resource is
always worse than assigning each a different resource in terms of
sum-throughput, we will focus on collision-free policies that assign
all users distinct resources, which we will refer to as a
permutation or matching. There are $P(N,M)$ such permutations.

We formulate our problem as a combinatorial multi-armed bandit, in
which each arm corresponds to a matching of the users to resources.
We can represent the arm corresponding to a permutation $k$ ($1 \leq
k \leq P(N,M)$) as the index set $\xA_k = \{(i,j): (i,j) \text{ is
in permutation } k \}$. The stochastic reward for choosing arm $k$
at time $n$ under policy $\alpha$ is then given as
\begin{equation}
Y_{\alpha(n)}(n) = \sum\limits_{(i,j) \in \xA_{\alpha(n)}}
Y_{i,j}(n) = \sum\limits_{(i,j) \in \xA_{\alpha(n)}}
\theta_{i,j}^{z_{\alpha(n)}}. \nonumber
\end{equation}
Note that different from most prior
work on multi-armed bandits, this combinatorial formulation results
in dependence across arms that share common components.

A key metric of interest in evaluating a given policy for this
problem is \emph{regret}, which is defined as the difference between
the expected reward that could be obtained by the best-possible
static matching, and that obtained by the given policy. It can be
expressed as:
\begin{equation}
\begin{split}
 R^\alpha (n) & = n \mu^*  - E_\alpha[ \sum \limits_{t = 1}^n Y_{\alpha(t)}(t) ]\\
 & = n \mu^*  - E_\alpha[ \sum \limits_{t = 1}^n \sum\limits_{(i,j) \in \xA_{\alpha(t)}} \theta_{i,j}^{z_{\alpha(t)}} ],
 \end{split}
\end{equation}
where $\mu^* = \max\limits_{k} \sum \limits_{(i,j) \in \xA_k}
\mu^{i,j}$, the expected reward of the optimal arm, is the expected
sum-weight of the maximum weight matching of users to resources with
$\mu_{i,j}$ as the weight.

We are interested in designing policies for this combinatorial
multi-armed bandit problem with Markovian rewards that perform well
with respect to regret. Intuitively, we would like the regret
$R^\alpha (n)$ to be as small as possible. If it is sub-linear with
respect to time $n$, the time-averaged regret will tend to zero.

\section{Matching Learning for Markovian Rewards} \label{sec:algorithm}

A straightforward idea for the combinatorial multi-armed bandit
problem with Markovian rewards is to treat each matching as an arm,
apply UCB1 policy (given by Auer~\emph{et al.}~\cite{Auer:2002})
directly, and ignore the dependencies across the different arms. For
each arm $k$, two variables are stored and updated: the time average
of all the observation values of arm $k$ and the number of times
that arm $k$ has been played up to the current time slot. The UCB1
policy makes decisions based on this information alone.

However, there are several problems that arise in applying UCB1
directly in the above setting. We note that UCB1 requires both the
storage and computation time that are linear in the number of arms.
Since the number of arms in this formulation grows as $P(N,M)$, it
is highly unsatisfactory. Also, the upper-bound of regret given in
\cite{Tekin:2010} will not work anymore since the rewards across
arms are not independent anymore and the states of an arm may
involve even when this arm is not played. No existing analytical
result on the upper-bound of regret can be applied directly in this
setting to the best of our knowledge.

So we are motivated to propose a policy which more efficiently
stores observations from correlated arms and exploits the
correlations to make better decisions. Our key idea is to use two
$M$ by $N$ matrices, $(\hat{ \theta }_{i,j})_{M \times N}$ and
$(n_{i,j})_{M \times N}$, to store the information for each
user-resource pair, rather than for each arm as a whole. $\hat{
\theta }_{i,j}$ is the average (sample mean) of all the observed
values of resource $j$ by user $i$ up to the current time slot
(obtained through potentially different sets of arms over time).
$n_{i,j}$ is the number of times that resource $j$ has been assigned
to user $i$ up to the current time slot.

At each time slot $n$, after an arm $k$ is played, we get the
observation of $Y_{i,j}(n)$ for all $(i,j) \in \xA_k$. Then $(\hat{
\theta }_{i,j})_{M \times N}$ and $(n_{i,j})_{M \times N}$ (both
initialized to 0 at time 0) are updated as follows:
\begin{equation} \label{eqn:update1}
 \hat{ \theta }_{i,j}(n) = \left\{ \begin{array}
  {l@{\;,\;}l}
  \frac{\hat{ \theta }_{i,j}(n-1) n_{i,j}(n-1) + Y_{i,j}(n)}{ n_{i,j}(n-1) +1} & \text{if } (i,j) \in \xA_k\\
  \hat{ \theta }_{i,j}(n-1) & \text{else} \\
  \end{array} \right.
\end{equation}

\begin{equation} \label{eqn:update2}
  n_{i,j}(n) = \left\{ \begin{array}
  {l@{\;,\;}l}
  n_{i,j}(n-1)+1 & \text{if } (i,j) \in \xA_k\\
  n_{i,j}(n-1) & \text{else} \\
  \end{array} \right.
\end{equation}

Note that while we indicate the time index in the above updates for
notational clarity, it is not necessary to store the matrices from
previous time steps while running the algorithm.

Our proposed policy, which we refer to as Matching Learning for
Markovian Rewards, is shown in Algorithm \ref{alg2:storage}.

\begin{algorithm} [ht]
\caption{Matching Learning for Markovian Rewards (MLMR)}
\label{alg2:storage}

\begin{algorithmic}[1]
\State \commnt{ Initialization}

\For {$p = 1$ to $M$}
    \For {$q = 1$ to $N$}
        \State $n = (M-1)p + q$;
        \State Play any permutation $k$ such that $(p,q) \in \xA_k$;
        \State Update $(\hat{\theta}_{i,j})_{M \times N}$, $(n_{i,j})_{M \times N}$ accordingly.
    \EndFor
\EndFor

\State \commnt{Main loop}

\While {1}
    \State $n = n + 1$;
    \State Solve the Maximum Weight Matching problem (e.g., using the Hungarian
    algorithm~\cite{Kuhn:1955}) on the bipartite graph of users and resources with
    edge weights $\left(\hat{ \theta }_{i,j} + \sqrt{ \frac{ L \ln n }{ n_{i,j}}
    }\right)_{M \times N}$ to play arm $k$ that maximizes
    \begin{equation}
    \label{equ:a2}
    \sum\limits_{(i,j) \in \xA_k} \left( \hat{ \theta }_{i,j} + \sqrt{ \frac{ L \ln n }{ n_{i,j}}
    } \right)
    \end{equation}
    where $L$ is a positive constant.
    \State Update $(\hat{\theta}_{i,j})_{M \times N}$, $(n_{i,j})_{M \times N}$ accordingly. 
\EndWhile
\end{algorithmic}
\end{algorithm}

\section{Analysis of Regret} \label{sec:regret}

We summarizes some notation we use in the description and analysis
of our MLMR policy in Table \ref{tab:notation}.

\begin{table}[htbp]
\centering \normalsize
\begin{tabular}{|l l|}
\hline
&\\
 $N$ : & number of resources.\\
 $M$ : & number of users, $M \leq N$.\\
 $k$ : & index of a parameter used for an arm, \\
 & $1 \leq k \leq P(N,M)$.\\
 $i,j$ : & index of a parameter used for user $i$, resource j.\\
 $\xA_k$ : & $\{(i,j): (i,j) \text{ is in permutation } k \}$\\
 $\mathcal{K}_{i,j}$: & $\{ \xA_k : (i,j) \in \xA_k  \}$\\
 $*$ : & index indicating that a parameter is for the \\
    & optimal arm. If there are multiple optimal arms, \\
    & $*$ refers to any of them. \\
 $n_{i,j}$: & number of times that resource $j$ has been\\
    & matched with user $i$ up to the current time slot.\\
 $\hat{ \theta }_{i,j}$: & average (sample mean) of all observed values\\
  & of resource $j$ by user $i$ up to current time slot.\\
 $n_i^k$: & $n_{i,j}$ such that $(i,j) \in \xA_k$ at current time slot.\\
 $S_{i,j}$: & state space of the Markov chain for \\ & user-resource pair $(i,j)$.\\
 $\mathbf{P}_{i,j}$: & transition matrix of the Markov chain \\
  & associated with user-resource pair $(i,j)$.\\
 $\pi^{i,j}_z$: & steady state distribution for state $z$ of the \\
  & Markov chain associated with $(i,j)$.\\
 $\theta^{i,j}_z$:  & reward obtained by user $i$ while access \\
    & resource $j$ on state $z \in S_{i,j}$. \\
 $ \mu_{i,j}$: & $\sum\limits_{z \in S^i} \theta^{i,j}_z
 \pi^{i,j}_z$, the mean reward for user $i$ using \\
    & resource $j$\\
 $\mu^{k}$: & $\sum\limits_{(i,j) \in \xA_k} \mu^{i,j}$\\
 $\mu^*$: & $\max\limits_{k} \sum \limits_{(i,j) \in \xA_k}
 \mu^{i,j}$\\
 $\Delta_k$: & $\mu^* - \mu_k$.\\
 $\Delta_{\min}$: & $\min\limits_{k: \mu_k < \mu^*} \Delta_k$.\\
 $\Delta_{\max}$: & $\max\limits_k \Delta_k$.\\
 $\pi_{\min}$: & $\min\limits_{1 \leq i
\leq M, 1 \leq j \leq N, z \in S_{i,j}} \pi^{i,j}_z$.\\
 $s_{\max}$: & $\max\limits_{1 \leq i \leq M, 1 \leq j \leq N} |S_{i,
 j}|$.\\
 $s_{\min}$: & $\min\limits_{1 \leq i \leq M, 1 \leq j \leq N} |S_{i,
j}|$.\\
 $\theta_{\max}$: & $\max\limits_{1 \leq i \leq M, 1 \leq j \leq
N, z \in S_{i,j}} \theta^{i,j}_z$.\\
 $\theta_{\min}$: & $ \min\limits_{1
 \leq i \leq M, 1 \leq j \leq N, z \in S_{i,j}} \theta^{i,j}_z$.\\
 $\epsilon_{i,j}$: & eigenvalue gap, defined as $1 - \lambda_2$, where $\lambda_2$ \\
 &  is the
 second largest eigenvalue of $\mathbf{P}_{i,j}$.\\
 $\epsilon_{\max}$: & $\max\limits_{1 \leq i \leq M, 1 \leq j \leq
 N} \epsilon_{i,j}$. \\
 $\epsilon_{\min}$: & $\min\limits_{1 \leq i \leq M, 1 \leq j \leq
 N} \epsilon_{i,j}$. \\
 $T_k (n)$: & number of times arm $k$ has been played by\\
 & MLMR in the first $n$ time slots.\\
 $\bX_k(n)$: & $\sum \limits_{(i,j) \in \xA_k} \hat{ \theta }_{i,j}(n)$. It is the summation of all the\\
 &  average observation values in arm $k$ at time $n$.\\
 $\hX_{i,n_i^k}^k$ : & $\hat{ \theta }_{i,j}(n)$ such that $(i,j) \in \xA_k$ and $n_{i,j}(n) = n_i^k$.\\
 \multicolumn{2}{|l|}{$\bX_{k, n_1^{k}, \ldots, n_M^{k}}$ : $\sum\limits_{i = 1}^M \hX_{i,n_i^k}^k$.}\\
\hline
\end{tabular}
\caption{Notation} \label{tab:notation}
\end{table}


The regret of a policy for a multi-armed bandit problem is
traditionally upper-bounded by analyzing the expected number of
times that each non-optimal arm is played and then taking the
summation over these expectation times the reward difference between
an optimal arm and a non-optimal arm all non-optimal arms. Although
we could use this approach to analyze the MLMR policy, we notice
that the upper-bound for regret consequently obtained is quite
loose, which is linear in the number of arms, $P(N,M)$. Instead, we
present here a novel analysis for a tighter analysis of the MLMR
policy. Our analysis shows an upper bound of the regret that is
polynomial in $M$ and $N$, and uniformly logarithmic over time.

The following lemmas are needed for our main results in Theorem
\ref{theorem:1}:

\lemma \label{lemma:1} (Lemma 2.1 from \cite{Anantharam:1987})
$\{X_n, n = 1, 2, \ldots\}$ is an irreducible aperiodic Markov chain
with state space $S$, transition matrix $P$, a stationary
distribution ${\pi_z}$, $\forall z \in S$, and an initial
distribution $\mathbf{q}$. Denote $F_t$ as the $\sigma$-algebra
generated by $X_1, X_2, \ldots, X_t$. Let $G$ be a $\sigma$-algebra
independent of $F = \vee_{t \geq 1} F_t$. Let $\tau$ be a stopping
time with respect to the increasing family of $\sigma$-algebra ${G
\vee F_t, t \geq 1}$. Define $N(z, \tau)$ such that $N(z, \tau) =
\sum\limits_{t = 1}^{\tau} I(X_t = z)$. Then,
\begin{equation}
 |E[N(z, \tau) - \pi_z E[ \tau ]]| \leq A_P,
\end{equation}
for all $\mathbf{q}$ and all $\tau$ such that $E[\tau ] < \infty$.
$A_P$ is a constant that depends on $P$.

\lemma \label{lemma:2} (Corollary 1 from \cite{Tekin:2010}) Let
$\pi_{\min}$ be the minimum value among the stationary distribution,
which is defined as $\pi_{\min} = \min\limits_{z \in S} \pi_z$. Then
$A_P \leq 1/\pi_{\min}$.

\lemma \label{lemma:3}

For user-resource matching, if the state of reward associated with
each user-resource pair $(i, j)$ is given by a Markov chain, denoted
$\{X^{i,j}_1, X^{i,j}_2, \ldots\}$, satisfying the properties of
Lemma \ref{lemma:1}, then the regret under policy $\alpha$ is
bounded by:
\begin{equation}
 R^\alpha (n) \leq \sum\limits_{k = 1}^{P(N,M)} (\mu^* - \mu^k) E_{\alpha}
 [{T_k^\alpha(n)}] + A_{\mathbf{S},\mathbf{P},\Theta},
\end{equation}
where $A_{\mathbf{S},\mathbf{P},\Theta}$ is a constant that depends
on all the state spaces $\{S_{i,j}\}_{1 \leq i \leq M, 1\leq i \leq
N}$, transition probability matrices $\{\mathbf{P}_{i,j}\}_{1 \leq i
\leq M, 1\leq i \leq N}$ and the rewards set $\{\theta^z_{i,j}, z
\in S_{i,j} \}_{1 \leq i \leq M, 1\leq i \leq N}$.

\begin{IEEEproof}

$\forall 1 \leq i \leq M$, $1 \leq j \leq N$, define $G_{i,j} =
\vee_{k \neq i, l \neq j} F_{k, l} $ where $F_{k, l} = \vee_{t \geq
1} F_t^{i,j}$, which applies to the Markove chain $\{X^{i,j}_1,
X^{i,j}_2, \ldots\}$. We note that the Markove chains of different
user-resource pairs are mutually independent, so $\forall i,j$,
$G^{i,j}$ is independent of $F_{i,j}$. $F_{i,j}$ satisfies the
conditions in Lemma \ref{lemma:1}. Note that $T_{i,j}^{\alpha}(n)$
is a stopping time with respect to $\{G_{i,j} \vee F_n^{i,j}, n
>1\}$.

Since the state of a Markove chain evolves only when it is observed,
$X^{i,j}_1, \ldots, X^{i,j}_{T_{i,j}^{\alpha}(n)}$ represents the
successive states of the Markov chain up to $n$ when assigning
resource $j$ to user $i$.Then the total reward obtained under policy
$\alpha$ up to time $n$ is given by:
\begin{equation}
\sum \limits_{t = 1}^n Y_{\alpha(t)}(t) = \sum\limits_{j = 1}^N
\sum\limits_{i = 1}^M \sum\limits_{l = 1}^{T_{i,j}^{\alpha}(n)}
\sum\limits_{z \in S_{i,j}} \theta^{i,j}_z I(X^{i,j}_l = z).
\end{equation}

Note that $\forall i = 1, \ldots, M$,
 $T_k^\alpha(n) = T_k^{\alpha(n),i}$ where $T_k^{\alpha(n),i}$ is the number of times up
to $n$ that the $i$-th component has been observed while playing arm
$k$, and there exist one resource index $j$ such that $(i, j) \in
\xA_k$. So, we have:

\begin{equation}
 \begin{split}
 & \sum\limits_{k = 1}^{P(N,M)} \mu^k E_{\alpha}
 [{T_k^\alpha(n)}]\\ \nonumber
 & = \sum\limits_{k = 1}^{P(N,M)} \sum\limits_{i =
 1}^M \mu_i^k E_{\alpha}
 [{T_k^\alpha(n)}] \\
 & = \sum\limits_{k = 1}^{P(N,M)} \sum\limits_{i =
 1}^M \mu_i^k E_{\alpha}
 [{T_k^{\alpha,i} (n)}]\\
 & = \sum\limits_{j
= 1}^N \sum\limits_{i = 1}^M \mu_{i,j} \sum\limits_{\xA_k \in
\mathcal{K}_{i,j} } E_{\alpha}
 [{T_k^{\alpha,i} (n)}] \\
 & = \sum\limits_{j
= 1}^N \sum\limits_{i = 1}^M \mu_{i,j} E_{\alpha}
 [{T_{i,j}^{\alpha} (n)}] \\
 & = \sum\limits_{j
= 1}^N \sum\limits_{i = 1}^M \sum\limits_{z \in S_{i,j}}
\theta^{i,j}_z \pi^{i,j}_z E_{\alpha}
 [{T_{i,j}^{\alpha} (n)}]
 \end{split}.
\end{equation}

Hence,
\begin{align}
 & |R^\alpha (n) - \sum\limits_{k = 1}^{P(N,M)} (\mu^* - \mu^k) E_{\alpha}
 [{T_k^\alpha(n)}]| \notag\\
 & = \left|R^\alpha (n) - ( n \mu^* - \sum\limits_{k = 1}^{P(N,M)} \mu^k E_{\alpha}
 [{T_k^\alpha(n)}] )\right| \notag\\
 & =  \left| ( n \mu^*  -E_\alpha[\sum \limits_{t = 1}^n Y_{\alpha(t)}(t) ]
 ) \right. \notag\\
 & \quad \left. - ( n \mu^* - \sum\limits_{k = 1}^{P(N,M)} \mu^k E_{\alpha}
 [{T_k^\alpha(n)}] )
 \right| \notag\\
 & = \left| E_\alpha[ \sum \limits_{t = 1}^n Y_{\alpha(t)}(t) ] - \sum\limits_{k = 1}^{P(N,M)} \mu^k E_{\alpha}
 [{T_k^\alpha(n)}] \right| \notag\\
 & = \left| E_\alpha[ \sum\limits_{j
= 1}^N \sum\limits_{i = 1}^M \sum\limits_{l =
1}^{T_{i,j}^{\alpha}(n)}
\sum\limits_{z \in S_{i,j}} \theta^{i,j}_z I(X^{i,j}_l = z) ] \right.\notag\\
& \quad \left. - \sum\limits_{j = 1}^N \sum\limits_{i = 1}^M
\sum\limits_{z \in S_{i,j}} \theta^{i,j}_z \pi^{i,j}_z E_{\alpha}
 [{T_{i,j}^{\alpha} (n)}]
\right| \notag\\
& \leq \sum\limits_{j = 1}^N \sum\limits_{i = 1}^M \sum\limits_{z
\in S_{i,j}} | E_\alpha[ \sum\limits_{l = 1}^{T_{i,j}^{\alpha}(n)}
\theta^{i,j}_z I(X^{i,j}_l = z) ] . \notag\\
& \quad  - \theta^{i,j}_z \pi^{i,j}_z E_{\alpha} [{T_{i,j}^{\alpha}
(n)}] | \notag\\
 & =  \sum\limits_{j = 1}^N \sum\limits_{i = 1}^M \sum\limits_{z
\in S_{i,j}} \theta^{i,j}_z | E_\alpha[ \sum\limits_{l =
1}^{T_{i,j}^{\alpha}(n)}
 I(X^{i,j}_l = z) ]  \notag\\
& \quad  - \pi^{i,j}_z E_{\alpha} [{T_{i,j}^{\alpha} (n)}] |\notag\\
 & =  \sum\limits_{j = 1}^N \sum\limits_{i = 1}^M \sum\limits_{z
\in S_{i,j}} \theta^{i,j}_z \left| E_{\alpha}[N(z, {T_{i,j}^{\alpha}
(n)})] -  \pi^{i,j}_z E_{\pi} [{T_{i,j}^{\alpha} (n)}] \right|.
\notag
\end{align}

Based on Lemma \ref{lemma:1}, we have:
\begin{align}
& |R^\alpha (n) - \sum\limits_{k = 1}^{P(N,M)} (\mu^* - \mu^k)
E_{\alpha} [{T_k^\alpha(n)}]| \notag\\
& \leq \sum\limits_{j = 1}^N \sum\limits_{i = 1}^M \sum\limits_{z
\in S_{i,j}} \theta^{i,j}_z C_{P_{i,j}} =
A_{\mathbf{S},\mathbf{P},\Theta}. \label{equ:11}
\end{align}
\end{IEEEproof}


\lemma \label{lemma:4} (Theorem 2.1 from \cite{Gillman:1998}) Let
$\{X_n, n = 1, 2, \ldots\}$ be an irreducible aperiodic Markov chain
with finite state space $S$, transition matrix $\mathbf{P}$, a
stationary distribution ${\pi_z}$, $\forall z \in S$, and an an
initial distribution $\mathbf{q}$. Let $N_{\mathbf{q}} = ||(
\frac{q_z}{\pi_z} ), z \in S||_2$. The eigenvalue gap $\epsilon$ is
defined as $\epsilon = 1 - \lambda_2$, where $\lambda_2$ is the
second largest eigenvalue of the matrix $\mathbf{P}$. $\forall A
\subset S$, define $t_A(n)$ as the total number of times that all
states in the set $A$ are visited up to time $n$. Then $\forall
\gamma \geq 0$,
\begin{equation}
 P(t_A(n) - n \pi_A \geq \gamma) \leq (1 + \frac{\gamma \epsilon}{10 n} N_{\mathbf{q}} e^{ -\gamma^2 \epsilon / 20 n
 }),
\end{equation}
where $\pi_A = \sum\limits_{z \in A} \pi_z$.

Our main results on the regret of MLMR policy are shown in Theorem
\ref{theorem:1}. We show that with using a constant $L$ which is
bigger than a value determined by the minimum eigenvalue gap of the
transition matrix, maximum value of the number of states, and
maximum value of the rewards, our MLMR policy is guaranteed to
achieve a regret that is uniformly logarithmic in time, and
polynomial in the number of users and resources.

\theorem\label{theorem:1} When using any constant $L \geq \frac{
(50+40M) \theta_{\max}^2 s^2_{\max} }{ \epsilon_{min} }$, the
expected regret under the MLMR policy specified in Algorithm
\ref{alg2:storage} is at most
\begin{equation}
\begin{split}
& \left[\frac{ 4 M^3 N L \ln n }{ \left( \Delta_{\min} \right)^2 } +
MN + \right.\\
& \left. M^2 N \frac{s_{\max}}{\pi_{\min}} \left( 1 +
\frac{\epsilon_{\max}\sqrt{L}}{10
 s_{\min} \theta_{\min}} \right) \frac{\pi}{3} \right] \Delta_{\max} +
 A_{\mathbf{S},\mathbf{P},\Theta},
 \end{split}
\end{equation}
where $\Delta_{\min}$, $\Delta_{\max}$, $\pi_{\min}$, $s_{\max}$,
$s_{\min}$, $\theta_{\max}$, $\theta_{\min}$, $\epsilon_{\max}$,
$\epsilon_{\min}$ follow the definition in Table \ref{tab:notation};
$A_{\mathbf{S},\mathbf{P},\Theta}$ follows the definition in Lemma
\ref{lemma:3}.


\begin{IEEEproof}

Denote $C_{t, n}$ as $\sqrt{ \frac{ L \ln t }{ n} }$. Denote $C_{t,
\mathbf{n}_{A_k}} = \sum\limits_{(i,j) \in \xA_k} \sqrt{ \frac{ L
\ln t }{ n_{i,j}} } = \sum\limits_{i = 1}^M \sqrt{ \frac{ L \ln t }{
n^k_i}}  = \sum\limits_{i = 1}^M C_{t, n^k_i}$. It is also denoted
as $C_{t, (n^k_1, \ldots, n^k_M)}$ sometimes for a clear explanation
in this proof.

We introduce $\tT(n)$ as a counter after the initialization period.
It is updated in the following way:

At each time slot after the initialization period, one of the two
cases must happen: (1) an optimal arm is played; (2) a non-optimal
arm is played. In the first case, $(\tT(n))_{M \times N}$ won't be
updated. When an non-optimal arm $k(n)$ is picked at time $n$, there
must be at least one $(i,j) \in \xA_k$ such that $n_{i,j}(n) =
\min\limits_{ (i_1,j_1) \in \xA_k } n_{i_1,j_1}$. If there is only
one such arm, $\tT(n)$ is increased by $1$. If there are multiple
such arms, we arbitrarily pick one, say $(i',j')$, and increment
$\widetilde{T}_{i'j'}$ by $1$.

Each time when a non-optimal arm is picked, exactly one element in
$(\tT(n))_{M \times N}$ is incremented by $1$. This implies that the
total number that we have played the non-optimal arms is equal to
the summation of all counters in $(\tT(n))_{M \times N}$. Therefore,
we have:
\begin{equation}
 \label{equ:f1}
 \sum\limits_{k: \mu_k < \mu*} \mathds{E}[T_k(n)] = \sum\limits_{i = 1}^M  \sum\limits_{j = 1}^N
 \mathds{E}[\tT(n)].
\end{equation}

Also note for $\tT(n)$, the following inequality holds:
\begin{equation}
 \label{equ:f2}
 \tT(n) \leq n_{i,j}(n), \forall 1 \leq i \leq M, 1 \leq j \leq N.
\end{equation}

Denote by $\tI(n)$ the indicator function which is equal to $1$ if
$\tT(n)$ is added by one at time $n$. Let $l$ be an arbitrary
positive integer. Then:
\begin{equation}
 \begin{split}
 \tT(n) & = \sum\limits_{t = MN+1}^n \mathds{1} \{ \tI(t)\} \nonumber \\
 & \leq l + \sum\limits_{t = MN+1}^n \mathds{1} \{ \tI(t) , \tT(t-1) \geq l \}
 \end{split}
\end{equation}
where $\mathds{1}(x)$ is the indicator function defined to be 1 when
the predicate $x$ is true, and 0 when it is false.


When $\tI(t) = 1$, there exists some arm such that a non-optimal arm
is picked for which $n_{i,j}$ is the minimum in this arm. We denote
this arm as $k(t)$ since at each time that $\tI(t) = 1$, we may get
different arms. Then,
\begin{equation}
 \begin{array}{r@{}l}
 \tT(n) & \leq l + \sum\limits_{t = MN+1}^n \mathds{1} \{ \obX(t-1) + C_{t-1, \mathbf{n}^*(t-1)}  \nonumber\\
  \leq & \bX_{k(t-1)}(t-1) + C_{t-1, \mathbf{n}_{\xA_{k(t-1)} }(t-1)}, \tT(t-1) \geq l \} \\
  = &\;  l + \sum\limits_{t = MN}^n \mathds{1} \{ \obX(t) + C_{t, \mathbf{n}^*(t)}  \\
  \leq  & \; \bX_{k(t)}(t) + C_{t, \mathbf{n}_{\xA_{k(t)} }(t)}, \tT(t) \geq l \}. \\
 \end{array}
\end{equation}
Based on (\ref{equ:f2}), $l \leq \tT(t)$ implies,
\begin{equation}
 l \leq \tT(t) \leq n_{i,j}(t) = n_i^{k(t)}. \nonumber
\end{equation}
So,
\begin{equation}
 \forall 1 \leq i \leq M, n_i^{k(t)} \geq l. \nonumber
\end{equation}
Then we could bound $\tT(n)$ as,
\begin{equation}
 \begin{array}{r@{}l}
 \tT(n) & \leq l + \sum\limits_{t = MN}^n \mathds{1} \{\min\limits_{0 < n_1^*, \ldots, n_M^* \leq t } \obX_{n_1^*, \ldots, n_M^*} \nonumber\\
 & + C_{t, (n_1^*, \ldots, n_M^*)} \leq \max\limits_{l \leq n_1^{k(t)}, \ldots, n_M^{k(t)} \leq t} \bX_{k(t), n_1^{k(t)}, \ldots, n_M^{k(t)}} \\
 & \quad + C_{t, (n_1^{k(t)}, \ldots, n_M^{k(t)})}\} \\[4mm]
 & \leq l + \sum\limits_{t = 1}^{\infty} [\sum\limits_{n_1^* = 1}^{t} \dots \sum\limits_{n_M^* = 1}^{t} \sum\limits_{n_1^{k(t)} = l}^{t} \dots \sum\limits_{n_M^{k(t)} = l}^{t}   \\
 & \quad \mathds{1}\{ \obX_{n_1^*, \ldots, n_M^*} + C_{t, (n_1^*, \ldots, n_M^*)} \leq \bX_{k(t), n_1^{k(t)}, \ldots, n_M^{k(t)}} \\
 & \quad + C_{t, (n_1^{k(t)}, \ldots, n_M^{k(t)})}\}]. \\
 \end{array}
\end{equation}

$\obX_{n_1^*, \ldots, n_M^*} + C_{t, (n_1^*, \ldots, n_M^*)} \leq
\bX_{k(t), n_1^{k(t)}, \ldots, n_M^{k(t)}} + C_{t, (n_1^{k(t)},
\ldots, n_M^{k(t)})} $ means that at least one of the following must
be true:
\begin{equation}
\label{equ:p1}
 \obX_{n_1^*, \ldots, n_M^*} \leq \mu^* - C_{t, (n_1^*, \ldots,
 n_M^*)},
\end{equation}
\begin{equation}
\label{equ:p2}
 \bX_{k(t), n_1^{k(t)}, \ldots, n_M^{k(t)}} \geq \mu_{k(t)} + C_{t, (n_1^{k(t)}, \ldots,
 n_M^{k(t)})},
\end{equation}
\begin{equation}
\label{equ:p3}
 \mu^* < \mu_{k(t)}  + 2 C_{t, (n_1^{k(t)}, \ldots, n_M^{k(t)})}.
\end{equation}

Here we first find the upper bound for $Pr\{ \obX_{n_1^*, \ldots,
n_M^*} \leq \mu^* - C_{t, (n_1^*, \ldots, n_M^*)} \}$:

\begin{equation}
 \begin{array}{r@{}l}
 & Pr\{ \obX_{n_1^*, \ldots, n_M^*} \leq \mu^* - C_{t, (n_1^*, \ldots, n_M^*)} \} \nonumber \\
 & = Pr\{ \sum\limits_{i = 1}^M \hX_{i,n_i^*}^* \leq \sum\limits_{i = 1}^M \mu_i^* - \sum\limits_{i = 1}^M  C_{t, n^*_i} \} \\
 & \leq Pr\{ \text{At least one of the following must hold:} \\
 & \qquad \qquad \hX_{1,n_1^*}^* \leq \mu_1^* - C_{t, n^*_1 }, \\
 & \qquad \qquad \hX_{2,n_2^*}^* \leq \mu_2^* - C_{t, n^*_2 }, \\
 & \qquad \qquad \qquad \qquad \vdots \\
 & \qquad \qquad \hX_{M,n_M^*}^* \leq \mu_M^* - C_{t, n^*_M } \} \\
 & \leq \sum\limits_{i = 1}^{M} Pr\{\hX_{i,n_i^*}^* \leq \mu_i^* -  C_{t, n^*_i }
 \}.
 \end{array}
\end{equation}
$\forall 1 \leq i \leq M$,
\begin{equation} \label{equ:b1}
 \begin{split}
 & Pr\{\hX_{i,n_i^*} \leq \mu_i^* - C_{t, n^*_i }
 \} \\
 & = Pr \{ \sum\limits_{z = 1}^{|S^*_i|} \frac{\theta_i^*(z) (z) n_i^*(z) }{ n_i^* }
 \leq  \sum\limits_{z = 1}^{|S^*_i|} \theta_i^* (z) \pi_i^* (z) -  C_{t, n^*_i }
 \} \\
 & = Pr \{ \sum\limits_{z = 1}^{|S^*_i|} ( \theta_i^*(z) n_i^*(z) - n_i^* \theta_i^*(z) \pi_i^* (z))
 \leq - n_i^* C_{t, n^*_i } \} \\
 & \leq Pr\{ \text{At least one of the following must hold:} \\
 &   \theta_i^*(1) n_i^*(1) - n_i^* \theta_i^*(1) \pi_i^*(1) \leq - \frac{n_i^*}{|S^*_i|} C_{t, n^*_i } , \\
 &  \qquad \qquad \vdots \\
&  \theta_i^* (|S^*_i|) n_i^*(|S^*_i|) - n_i^* \theta_i^*
 (|S^*_i|) \pi_i^*(|S^*_i|) \leq - \frac{n_i^*}{|S^*_i|} C_{t, n^*_i } \}, \\
 & \leq  \sum\limits_{z = 1}^{|S^*_i|} Pr\{ \theta_i^*(z) n_i^*(z) - n_i^* \theta_i^*(z)
  \pi_i^*(z) \leq - \frac{n_i^*}{|S^*_i| } C_{t, n^*_i } \}\\
 & = \sum\limits_{z = 1}^{|S^*_i|} Pr\{  n_i^*(z) - n_i^*
  \pi_i^*(z) \leq - \frac{n_i^*}{|S^*_i| \theta_i^*(z)} C_{t, n^*_i } \} \\
 & = \sum\limits_{z = 1}^{|S^*_i|} Pr\{ (n_i^*- \sum\limits_{l \neq z }n_i^*(l)) -
 n_i^* (1- \sum\limits_{l \neq z } \pi_i^*(z)) \\
 & \qquad \leq - \frac{n_i^*}{|S^*_i| \theta_i^*(z)} C_{t,  n^*_i } \} \\
 & = \sum\limits_{z = 1}^{|S^*_i|} Pr\{ \sum\limits_{l \neq z }n_i^*(l) - n_i^* \sum\limits_{l \neq z } \pi_i^*(z)
 \geq \frac{n_i^*}{|S^*_i| \theta_i^*(z)} C_{t, n^*_i } \}.
 \\
  \end{split}
\end{equation}

$\forall 1 \leq z \leq |S^*_i|$, applying Lemma \ref{lemma:4}, we
could find the upper bound of each probablilty in (\ref{equ:b1}) as,

\begin{align}
 & Pr\{\hX_{i,n_i^*} \leq \mu_i^* - C_{t, n^*_i }
 \} \notag\\
 & \leq \sum\limits_{z = 1}^{|S^*_i|} \left( 1 + \frac{\epsilon_{i,j}}{10 |S^*_i|
 \theta_i^*(z)} \right)
 \sqrt{ \frac{L \ln t}{n^*_i} }  N_{\mathbf{q}_{i,j}} e^{
 -\frac{ n_i^* L \ln t \epsilon_{i,j} }{ 20 |S^*_i|^2 \theta_i^*(z)^2  n^*_i}}
  \notag\\
 & \leq \sum\limits_{z = 1}^{|S^*_i|} \left( 1 + \frac{\epsilon_{\max}\sqrt{L t}}{10
 s_{\min} \theta_{\min}} \right) N_{\mathbf{q}_{i,j}} e^{
 -\frac{ L \ln t \epsilon_{\min} }{ 20 s_{\max}^2 \theta_{\max}^2 }}
  \notag\\
 & \leq \frac{s_{\max}}{\pi_{\min}} \sqrt{t} \left( 1 + \frac{\epsilon_{\max}\sqrt{L}}{10
 s_{\min} \theta_{\min}} \right) t^{-\frac{L \epsilon_{\min} }{20 s_{\max}^2 \theta_{\max}^2 } } \label{equ:20}\\
 & = \frac{s_{\max}}{\pi_{\min}} \left( 1 + \frac{\epsilon_{\max}\sqrt{L}}{10
 s_{\min} \theta_{\min}} \right) t^{-\frac{L \epsilon_{\min} - 10 s_{\max}^2 \theta_{\max}^2  }{20 s_{\max}^2 \theta_{\max}^2 }
 }, \notag
\end{align}
where (\ref{equ:20}) holds since for any $\mathbf{q}_{i,j}$,
\begin{equation}
\begin{split}
 N_{\mathbf{q}_{i,j}} & = \left\|\frac{q^{i,j}_z}{\pi^{i,j}_z}, z \in S_{i,j} \right\|_2 \leq \sum\limits_{z = 1}^{|S_{i,j}|}\left\|\frac{q^{i,j}_z}{\pi^{i,j}_z}
 \right\|_2\\
 & \leq \sum\limits_{z = 1}^{|S_{i,j}|}\frac{\left\|
 q^{i,j}_z\right\|_2}{\pi_{\min}} = \frac{1}{\pi_{\min}}. \nonumber
\end{split}
\end{equation}

Thus,
\begin{equation}\label{equ:23}
\begin{split}
 & Pr\{ \obX_{n_1^*, \ldots, n_M^*} \leq \theta^* - C_{t, (n_1^*, \ldots, n_M^*)}
 \}\\
 & \leq \frac{M s_{\max}}{\pi_{\min}} \left( 1 + \frac{\epsilon_{\max}\sqrt{L}}{10
 s_{\min} \theta_{\min}} \right) t^{-\frac{L \epsilon_{\min} - 10 s_{\max}^2 \theta_{\max}^2  }{20 s_{\max}^2 \theta_{\max}^2 } }.
\end{split}
\end{equation}

With the similar calculation, we can also get the upper bound of the
probability for (\ref{equ:p2}):
\begin{equation}\label{equ:24}
\begin{split}
& Pr \{  \bX_{k(t), n_1^{k(t)}, \ldots, n_M^{k(t)}} \geq \mu_k +
C_{t, (n_1^{k(t)}, \ldots, n_M^{k(t)})} \} \\
& \leq \sum\limits_{i = 1}^{M}  Pr\{\hX_{i,n_i^k}^k \geq \mu_i^k +
 C_{t, n^k_i } \}\\
& = \sum\limits_{i = 1}^{M} Pr \{ \sum\limits_{z = 1}^{|S^k_i|}
\frac{\theta_i^k(z) n_i^k(z) }{ n_i^k }
 \geq  \sum\limits_{z = 1}^{|S^k_i|} \theta_i^k (z) \pi_i^k (z) +  C_{t, n^k_i }
 \} \\
& \leq \sum\limits_{i = 1}^{M} \sum\limits_{z = 1}^{|S^k_i|} Pr\{
\theta_i^k(z) n_i^k(z) - n_i^k \theta_i^k(z)
  \pi_i^k(z) \geq \frac{n_i^k}{|S^*_i|} C_{t, n^k_i } \} \\
& = \sum\limits_{i = 1}^{M} \sum\limits_{z = 1}^{|S^k_i|} Pr\{
n_i^k(z) - n_i^k
  \pi_i^k(z) \geq \frac{n_i^k}{|S^k_i| \theta_i^k(z)} C_{t, n^k_i } \} \\
 & \leq \sum\limits_{i = 1}^{M} \frac{s_{\max}}{\pi_{\min}} \left( 1
 + \frac{\epsilon_{\max}\sqrt{L}}{10
 s_{\min} \theta_{\min}} \right) t^{-\frac{L \epsilon_{\min} - 10 s_{\max}^2 \theta_{\max}^2  }{20 s_{\max}^2 \theta_{\max}^2 } } \\
 & \leq \frac{M s_{\max}}{\pi_{\min}} \left( 1 +
 \frac{\epsilon_{\max}\sqrt{L}}{10
 s_{\min} \theta_{\min}} \right) t^{-\frac{L \epsilon_{\min} - 10 s_{\max}^2 \theta_{\max}^2  }{20 s_{\max}^2 \theta_{\max}^2 } }.
\end{split}
\end{equation}

Note that for $l \geq \left\lceil \frac{4 L \ln n}{
\left(\frac{\Delta_{k(t)}}{M} \right)^2 } \right\rceil$,
\begin{equation}
 \label{equ:p4}
\begin{split}
 & \mu^* - \mu_{k(t)} - 2 C_{t, (n_1^{k(t)}, \ldots, n_M^{k(t)})}  \\
  & =  \mu^* - \mu_{k(t)} - 2 \sum\limits_{i = 1}^M \sqrt{ \frac{ L \ln t }{ n^{k(t)}_i}}  \\
  & \geq  \mu^* - \mu_{k(t)} - M \sqrt{ \frac{ 4 L \ln n }{ 4 L \ln n } \left( \frac{\Delta_{k(t)} }{M} \right)^2 } \\
  &  =  \mu^* - \mu_{k(t)} - \Delta_{k(t)} = 0.
 \end{split}
\end{equation}
(\ref{equ:p4}) implies that condition (\ref{equ:p3}) is false when
$l = \left\lceil \frac{4L \ln n}{ \left(\frac{\Delta_{k(t)}}{M}
\right)^2 } \right\rceil$. If we let $l = \left\lceil \frac{4L \ln
n}{ \left(\frac{ \Delta_{\min}^{i,j} }{M} \right)^2 } \right\rceil$,
then (\ref{equ:p3}) is false for all $k(t), 1 \leq t \leq \infty$
where
\begin{equation}\label{equ:21}
\Delta_{\min}^{i,j} = \min\limits_k \{ \Delta_k: (i,j) \in  \xA_k
\}.
\end{equation}

Therefore,
\begin{align}
 & \mathds{E}[\tT(n)] \notag\\
 & \leq \left\lceil \frac{ 4L \ln n }{ \left( \frac{\Delta_{\min}^{i,j}}{M} \right)^2 } \right\rceil
 + \sum\limits_{t = 1}^\infty \left(\sum\limits_{n_1^* = 1}^{t} \dots \sum\limits_{n_1^* = M}^{t} \sum\limits_{n_1^k = 1}^{t} \dots \sum\limits_{n_1^k = M}^{t} \right.\notag\\
 & \left. 2 M
 \frac{s_{\max}}{\pi_{\min}} \left( 1 + \frac{\epsilon_{\max}\sqrt{L}}{10
 s_{\min} \theta_{\min}} \right) t^{-\frac{L \epsilon_{\min} - 10 s_{\max}^2 \theta_{\max}^2  }{20 s_{\max}^2 \theta_{\max}^2 } } \right) \notag\\
 & \leq  \frac{ 4 M^2 L \ln n }{ \left( \Delta_{\min}^{i,j} \right)^2 } + 1 \notag\\
 & + M \frac{s_{\max}}{\pi_{\min}} \left( 1 + \frac{\epsilon_{\max}\sqrt{L}}{10
 s_{\min} \theta_{\min}} \right) \sum\limits_{t = 1}^\infty 2 t^{-\frac{L \epsilon_{\min} - (40 M + 10) s_{\max}^2 \theta_{\max}^2  }{20 s_{\max}^2 \theta_{\max}^2 } } \notag\\
 & \leq  \frac{ 4 M^2 L \ln n }{ \left( \Delta_{\min}^{i,j} \right)^2 } + 1 + M \frac{s_{\max}}{\pi_{\min}} \left( 1 + \frac{\epsilon_{\max}\sqrt{L}}{10
 s_{\min} \theta_{\min}} \right) \sum\limits_{t = 1}^\infty 2 t^{-2} \label{equ:25}\\
 & = \frac{ 4 M^2 L \ln n }{ \left( \Delta_{\min}^{i,j} \right)^2 } + 1 + M \frac{s_{\max}}{\pi_{\min}} \left( 1 + \frac{\epsilon_{\max}\sqrt{L}}{10
 s_{\min} \theta_{\min}} \right) \frac{\pi}{3}, \notag
\end{align}
where (\ref{equ:25}) holds since $L \geq \frac{ (50+40M)
\theta_{\max}^2 s^2_{\max} }{ \epsilon_{min} }$.

So under our MLMR policy,
\begin{equation}
 \label{equ:p6}
\begin{split}
 & R_\pi (n) \leq \sum\limits_{k = 1}^{P(N,M)} (\mu^* - \mu^k) E_{\pi}
 [{T^k_\pi(n)}] + A_{\mathbf{S},\mathbf{P},\Theta} \\
 & = \sum\limits_{k: \theta_k < \theta*} \Delta_k \mathds{E}[T_k(n)] + A_{\mathbf{S},\mathbf{P},\Theta}\\
 & \leq \Delta_{\max} \sum\limits_{k: \theta_k < \theta*} \mathds{E}[T_k(n)] + A_{\mathbf{S},\mathbf{P},\Theta}\\
 & = \Delta_{\max} \sum\limits_{i = 1}^M  \sum\limits_{j = 1}^N \mathds{E}[\tT(n)] + A_{\mathbf{S},\mathbf{P},\Theta}\\
  \end{split}
\end{equation}
\begin{equation}
 \label{equ:p62}
\begin{split}
 & \leq \left[\sum\limits_{i = 1}^M  \sum\limits_{j = 1}^N  \frac{ 4 M^2 L \ln n }{ \left( \Delta_{\min}^{i,j} \right)^2 } + 1 \right. \\
   & \left. + M \frac{s_{\max}}{\pi_{\min}} \left( 1 + \frac{\epsilon_{\max}\sqrt{L}}{10
 s_{\min} \theta_{\min}} \right) \frac{\pi}{3} \right] \Delta_{\max} + A_{\mathbf{S},\mathbf{P},\Theta}\\
 & \leq \left[\frac{ 4 M^3 N L \ln n }{ \left( \Delta_{\min} \right)^2 } + M N \right.\\
   & \left. + M^2 N \frac{s_{\max}}{\pi_{\min}} \left( 1 + \frac{\epsilon_{\max}\sqrt{L}}{10
 s_{\min} \theta_{\min}} \right) \frac{\pi}{3} \right] \Delta_{\max} + A_{\mathbf{S},\mathbf{P},\Theta}.\\
  \end{split}
\end{equation}

\end{IEEEproof}

Theorem \ref{theorem:1} shows when we use a constant $L$ which is
large enough such that $L \geq \frac{ (50+40M) \theta_{\max}^2
s^2_{\max} }{ \epsilon_{min} }$, the regret of Algorithm
\ref{alg2:storage} is upper-bounded uniformly over time $n$ by a
function that grows as $O(M^3 N \ln n)$. However, when
$\theta_{\max}$, $s_{\max}$ or $\epsilon_{min}$ is unknown, the
upper bound of regret could not be guaranteed to grow
logarithmically in $n$.

So when no knowledge about the system is available, we extend the
MLMR policy to achieve a regret that is bounded uniformly over time
$n$ by a function that grows as $O(M^3 N L(n) \ln n)$, by using any
arbitrarily slowly diverging non-decreasing sequence $L(n)$ such
that $L(n) \leq n $ for any $n$ in Algorithm \ref{alg2:storage}.
Since $L(n)$ could grow arbitrarily slowly, the MLMR could achieve a
regret arbitrarily close to the logarithmic order. We present our
analysis in Theorem \ref{theorem:2}.

\theorem\label{theorem:2} When using any arbitrarily slowly
diverging non-decreasing sequence $L(n)$ (i.e., $L(n)\rightarrow
\infty$ as $n \rightarrow \infty$) in (\ref{equ:a2}) such that
$\forall n$, $L(n) \leq n $, the expected regret under the MLMR
policy specified in Algorithm \ref{alg2:storage} is at most
\begin{equation}
\begin{split}
& \left[\frac{ 4 M^3 N L(n) \ln n }{ \left( \Delta_{\min} \right)^2 } + M N B_{\mathbf{S},\mathbf{P},\Theta} \right.\\
   & \left. + M^2 N \frac{s_{\max}}{\pi_{\min}} \left( 1 + \frac{\epsilon_{\max}}{10
 s_{\min} \theta_{\min}} \right) \frac{\pi}{3} \right] \Delta_{\max} + A_{\mathbf{S},\mathbf{P},\Theta},
 \end{split}
\end{equation}
where $B_{\mathbf{S},\mathbf{P},\Theta}$ is a constant that depends
on $\theta_{\max}$, $s_{\max}$ and $\epsilon_{min}$.

\begin{IEEEproof}

Denote $C_{t, n}$ as $\sqrt{ \frac{ L(t) \ln t }{ n} }$. Denote
$C_{t, \mathbf{n}_{A_k}}$ as $ \sum\limits_{(i,j) \in \xA_k} \sqrt{
\frac{ L(t) \ln t }{ n_{i,j}} }$. Then replacing $L$ with $L(t)$ in
the proof of Theorem \ref{theorem:1}, (\ref{equ:f1}) to
(\ref{equ:21}) still stand. The upper bound of $\mathds{E}[\tT(n)]$
in (\ref{equ:25}) should be modified as in (\ref{equ:26}).

$L(t)$ is a diverging non-decreasing sequence, so there exists a
constant $t_1$, such that for all $t \geq t_1$, $L(t) \geq \frac{
(60+40M) \theta_{\max}^2 s^2_{\max} }{ \epsilon_{min} }$, which
implies $t^{-\frac{L(t) \epsilon_{\min} - (40 M + 10) s_{\max}^2
\theta_{\max}^2  }{20 s_{\max}^2 \theta_{\max}^2 } + \frac{1}{2}}
\leq t^{-2}$.

\begin{figure*}[!t]
\normalsize
\setcounter{MYtempeqncnt}{\value{equation}}
\setcounter{equation}{25}
\begin{align}
  \mathds{E}[\tT(n)]
  & \leq \left\lceil \frac{ 4L(n) \ln n }{ \left( \frac{\Delta_{\min}^{i,j}}{M} \right)^2 } \right\rceil
  + \sum\limits_{t = 1}^\infty \left(\sum\limits_{n_1^* = 1}^{t} \dots \sum\limits_{n_1^* = M}^{t} \sum\limits_{n_1^k = 1}^{t} \dots \sum\limits_{n_1^k = M}^{t}
   2 M \frac{s_{\max}}{\pi_{\min}} \left( 1 + \frac{\epsilon_{\max}\sqrt{L(t)}}{10
 s_{\min} \theta_{\min}} \right) t^{-\frac{L(t) \epsilon_{\min} - 10 s_{\max}^2 \theta_{\max}^2  }{20 s_{\max}^2 \theta_{\max}^2 } } \right) \notag\\
 & \leq  \frac{ 4 M^2 L(n) \ln n }{ \left( \Delta_{\min}^{i,j} \right)^2 } + 1 +
   M \frac{s_{\max}}{\pi_{\min}} \left( 1 + \frac{\epsilon_{\max}}{10
 s_{\min} \theta_{\min}} \right) \sum\limits_{t = 1}^\infty 2 \sqrt{L(t)} t^{-\frac{L(t) \epsilon_{\min} - (40 M + 10) s_{\max}^2 \theta_{\max}^2  }{20 s_{\max}^2 \theta_{\max}^2 } } \notag\\
 & \leq  \frac{ 4 M^2 L(n) \ln n }{ \left( \Delta_{\min}^{i,j} \right)^2 } + 1 +
   M \frac{s_{\max}}{\pi_{\min}} \left( 1 + \frac{\epsilon_{\max}}{10
 s_{\min} \theta_{\min}} \right) \sum\limits_{t = 1}^\infty 2 t^{-\frac{L(t) \epsilon_{\min} - (40 M + 10) s_{\max}^2 \theta_{\max}^2  }{20 s_{\max}^2 \theta_{\max}^2 } + \frac{1}{2}} \label{equ:26}
\end{align}
\setcounter{equation}{26}
\hrulefill
\vspace*{4pt}
\end{figure*}

Thus, we have

\begin{equation}
\begin{split}
 & \mathds{E}[\tT(n)] \leq  \frac{ 4 M^2 L(n) \ln n }{ \left( \Delta_{\min}^{i,j} \right)^2 }\\
 & \quad + M \frac{s_{\max}}{\pi_{\min}} \left( 1 + \frac{\epsilon_{\max}}{10
 s_{\min} \theta_{\min}} \right) \sum\limits_{t = t_1}^\infty 2 t^{-2} + B_{\mathbf{S},\mathbf{P},\Theta}\\
 & = \frac{ 4 M^2 L(n) \ln n }{ \left( \Delta_{\min}^{i,j} \right)^2 } + M \frac{s_{\max}}{\pi_{\min}} \left( 1 + \frac{\epsilon_{\max}\sqrt{L}}{10
 s_{\min} \theta_{\min}} \right) \frac{\pi}{3} + B_{\mathbf{S},\mathbf{P},\Theta}
 \end{split}
\end{equation}
where $B_{\mathbf{S},\mathbf{P},\Theta}$ is a constant as shown in
(\ref{equ:27}), which depends on $\theta_{\max}$, $s_{\max}$ and
$\epsilon_{min}$.

\begin{figure*}[!t]
\normalsize

\setcounter{MYtempeqncnt}{\value{equation}}

\setcounter{equation}{27}

\begin{equation} \label{equ:27}
\begin{split}
 B_{\mathbf{S},\mathbf{P},\Theta}
 =  1 + M \frac{s_{\max}}{\pi_{\min}} \left( 1 +
\frac{\epsilon_{\max}}{10
 s_{\min} \theta_{\min}} \right) \sum\limits_{t = 1}^{t_1 -
1} 2 t^{-\frac{L(t) \epsilon_{\min} - (40 M + 10) s_{\max}^2
\theta_{\max}^2  }{20 s_{\max}^2 \theta_{\max}^2 } + \frac{1}{2}}.
 \end{split}
\end{equation}

\setcounter{equation}{28}
\hrulefill
\vspace*{4pt}
\end{figure*}

Then for the MLMR policy with $L(n)$,
\begin{equation}
 \label{equ:p6}
\begin{split}
 & R_\pi (n) \leq \Delta_{\max} \sum\limits_{i = 1}^M  \sum\limits_{j = 1}^N \mathds{E}[\tT(n)] + A_{\mathbf{S},\mathbf{P},\Theta}\\ \\
 & \leq \left[\frac{ 4 M^3 N L(n) \ln n }{ \left( \Delta_{\min} \right)^2 } + M N B_{\mathbf{S},\mathbf{P},\Theta} \right.\\
   & \left. + M^2 N \frac{s_{\max}}{\pi_{\min}} \left( 1 + \frac{\epsilon_{\max}}{10
 s_{\min} \theta_{\min}} \right) \frac{\pi}{3} \right] \Delta_{\max} + A_{\mathbf{S},\mathbf{P},\Theta}.\\
  \end{split}
\end{equation}

\end{IEEEproof}

\section{Examples and Simulation Results} \label{sec:simulation}


We consider a system that consists of $M = 2$ users and $N = 4$
resources. The state of each resource evolves as an irreducible,
aperiodic Markov chain with two states ``0'' and ``1''. For all the
tables in this section, the element in the $i$-th row and $j$-th
column represents the value for the user-resource pair $(i,j)$. The
transition probabilities are shown in the tables below:

\vspace{0.25cm}
\begin{tabular}{|c|c|c|c|}
 \hline
 0.5&0.4&0.7&0.3\\
 \hline
 0.2&0.9&0.9&0.7\\
 \hline
\multicolumn{4}{c}{$\textbf{p}_{01}$}
\end{tabular}
\qquad
\begin{tabular}{|c|c|c|c|}
 \hline
 0.6&0.7&0.8&0.9\\
 \hline
 0.9&0.5&0.4&0.4\\
 \hline
\multicolumn{4}{c}{$\textbf{p}_{10}$}
\end{tabular}

\vspace{0.08cm}

The rewards on each states are:

\vspace{0.25cm}
\begin{tabular}{|c|c|c|c|}
 \hline
 0.6&0.5&0.2&0.4\\
 \hline
 0.3&0.7&0.8&0.3\\
 \hline
\multicolumn{4}{c}{$\bm{\theta}_{0}$}
\end{tabular}
\qquad
\begin{tabular}{|c|c|c|c|}
 \hline
 0.8&0.2&0.7&0.5\\
 \hline
 0.5&0.3&0.6&0.6\\
 \hline
\multicolumn{4}{c}{$\bm{\theta}_{1}$}
\end{tabular}

\vspace{0.08cm}

For $1 \leq i \leq M$, $1 \leq j \leq N$, the stationary
distribution of user-resource pair $(i,j)$ on state ``0'' is
calculated as $\frac{p_{10}^{i,j}}{p_{01}^{i,j}+p_{10}^{i,j}}$;  the
stationary distribution on state ``1'' is calculated as
$\frac{p_{01}^{i,j}}{p_{01}^{i,j}+p_{10}^{i,j}}$. The eigenvalue gap
is $\epsilon_{i,j} = p_{01}^{i,j}+p_{10}^{i,j}$. The expected reward
$\mu_{i,j}$ for all the pairs can be calculated as:

\vspace{0.25cm} \qquad\quad
\begin{tabular}{|c|c|c|c|}
 \hline
 0.6909&0.3909&0.4333&0.425\\
 \hline
 0.3363&0.4429&0.6615&0.4909\\
 \hline
\multicolumn{4}{c}{$\bm{\mu}$}
\end{tabular}

\vspace{0.08cm}

We can see that the arm $\{(1,1),(2,3)\}$ is the optimal arm with
greatest expected reward $\mu^* = 0.6909+0.6615 = 1.3524$.
$\Delta_{\min} = 0.1706$.
\begin{figure}[ht]
  \centering
  \includegraphics[width=0.5\textwidth]{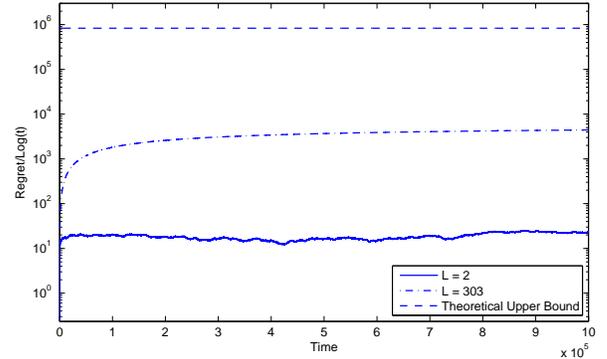}
 \caption{Simulation Results of Example 1 with $\Delta_{\min} = 0.1706$} \label{fig:1}
\end{figure}

Figure \ref{fig:1} shows the simulation result of the regret
(normalized with respect to the logarithm of time) for our MLMR
policy for the above system with different choices of $L$. We also
show the theoretical upper bound for comparison. The value of $L$ to
satisfy the condition in Theorem \ref{theorem:1} is $L \geq \frac{
(50+40M) R^2 s^2_{\max} }{ \epsilon_{min} } = 303$, so we picked $L
= 303$ in the simulation.

Note that in the proof of Theorem \ref{theorem:1}, when $L < \frac{
(50+40M) R^2 s^2_{\max} }{ \epsilon_{min} }$,
$-\frac{L \epsilon_{\min} - (40 M + 10) s_{\max}^2 \theta_{\max}^2
}{20 s_{\max}^2 \theta_{\max}^2 } \nonumber
> -2$.
This implies $\sum\limits_{t = 1}^\infty 2 t^{-\frac{L
\epsilon_{\min} - (40 M + 10) s_{\max}^2 \theta_{\max}^2  }{20
s_{\max}^2 \theta_{\max}^2
 }}$ does not converge anymore and thus we could not bound
 $\mathds{E}[\tT(n)]$ any more. Empirically, however, in \ref{fig:1}
 the case when $L=2$ also seems to yield logarithmic regret over
 time and the performance is in fact better
 than $L = 303$, since the non-optimal arms
 are played less when $L$ is smaller. However, this may possibly be due to
 the fact that the cases when $\tT(n)$ grows faster than $\ln(t)$
 only happens with very small probability when $L = 2$.

Table II shows the number of times that resource $j$ has been
matched with user $i$ up to time $n = 10^7$.

\vspace{0.25cm} \qquad\quad\;
\begin{tabular}{|c|c|c|c|}
 \hline
 999470&153&185&196\\
 \hline
 136&293&999155&420\\
 \hline
 \multicolumn{4}{c}{$n_{i,j}(10^7)$, $L = 2$}\\
\end{tabular}
\vspace{0.08cm}

\qquad\; \begin{tabular}{|c|c|c|c|}
 \hline
 892477&30685&39410&37432\\
 \hline
 26813&50341&850265&72585\\
 \hline
 \multicolumn{4}{c}{$n_{i,j}(10^7)$, $L = 303$}\\
 \multicolumn{4}{c}{\footnotesize{TABLE II} }
\end{tabular}



\begin{figure}[ht]
  \centering
  \includegraphics[width=0.5\textwidth]{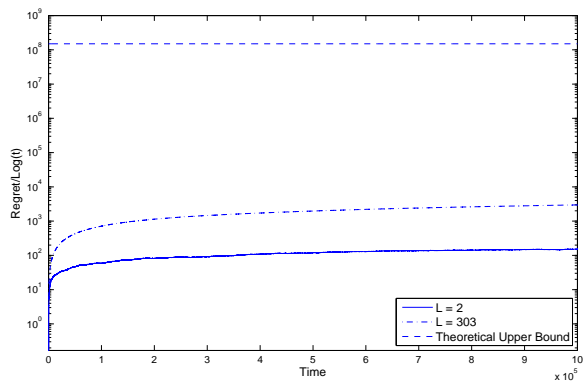}
 \caption{Simulation Results of Example 2 with $\Delta_{\min} = 0.0091$} \label{fig:2}
\end{figure}

Figure \ref{fig:2} shows the simulation results of the regret of
another example with the same transition probabilities as in the
previous example and different rewards on states as below:

\vspace{0.25cm}

\begin{tabular}{|c|c|c|c|}
 \hline
 0.7&0.3&0.5&0.5\\
 \hline
 0.65&0.7&0.8&0.4\\
 \hline
\multicolumn{4}{c}{$\bm{\theta}_{0}$}
\end{tabular}
\quad
\begin{tabular}{|c|c|c|c|}
 \hline
 0.4&0.6&0.7&0.45\\
 \hline
 0.5&0.5&0.6&0.55\\
 \hline
\multicolumn{4}{c}{$\bm{\theta}_{1}$}
\end{tabular}

\vspace{0.08cm}

The expected reward $\mu_{i,j}$ for all the pairs can be calculated
as:

\vspace{0.25cm} \qquad\quad
\begin{tabular}{|c|c|c|c|}
 \hline
 0.5636&0.4091&0.5933&0.4875\\
 \hline
 0.6227&0.5714&0.6615&0.4954\\
 \hline
 \multicolumn{4}{c}{$\bm{\mu}$}
\end{tabular}
\vspace{0.08cm}

$\{(1,1),(2,3)\}$ is still the optimal arm. However, compared with
the previous example, we can see that the expected reward of three
other arms $\{(1,3),(2,1)\}$, $\{(1,3),(2,2)\}$, $\{(1,1),(2,2)\}$
are all very close to the expected reward of the optimal arm. For
this example, $\Delta_{\min} = 0.0091$, which is much smaller
compared with the previous example. In this case, the non-optimal
arms are played much more compared with the previous example. This
is because we have several arms of which the expected rewards are
very close to $\mu^*$, so the policy has to spend a lot more time to
explore on those non-optimal arms to make sure those are non-optimal
arms. This fact can be seen clearly in Table III, which presents the
number of times that resource $j$ has been matched with user $i$ up
to time $n = 10^7$ under both cases when $L = 2$ and $L = 303$.

\vspace{0.25cm} \qquad\;\,
\begin{tabular}{|c|c|c|c|}
 \hline
 817529&544&179832&2099\\
 \hline
 175583&3610&820097&714\\
 \hline
 \multicolumn{4}{c}{$n_{i,j}(10^7)$, $L = 2$}\\
\end{tabular}
\vspace{0.08cm}

\qquad \begin{tabular}{|c|c|c|c|}
 \hline
 346395&60031&472346&121232\\
 \hline
 301491&146317&482545&69651\\
 \hline
 \multicolumn{4}{c}{$n_{i,j}(10^7)$, $L = 303$}\\
 \multicolumn{4}{c}{\footnotesize{TABLE III } }
\end{tabular}
\vspace{0.08cm}



\section{Conclusion}\label{sec:conclusion}

We have presented the MLMR policy for the problem of learning
combinatorial matchings of users to resources when the reward
process is Markovian. We showed that this policy requires only
polynomial storage and computation per step, and yields a regret
that grows uniformly logarithmically over time and only polynomially
with the number of users and resources.

In future work, we would like to also consider the case when the
rewards evolve not just when a user-resource pair is selected, but
rather at each discrete time. Further, we would like to investigate
if it is possible to analyze regret with respect to the best
non-static policy, which would be a stronger notion of regret than
that considered in this paper but is much harder to analyze.
Finally, exploring distributed schemes is also of interest, though
likely to be highly challenging in case of limited information
exchange between users.


\begin{thebibliography}{1}

\bibitem{Pandey}
S.~Pandey, D.~Chakrabarti, and D.~Agarwal, ``Multi-armed bandit
problems with dependent arms,'' In \emph{Proc. of the 24th
International Conference on Machine Learning }, Corvalis, June 2007.

\bibitem{Gai:2010}
Y.~Gai, B.~Krishnamachari, and R.~Jain, ``Learning multiuser channel
allocations in cognitive radio networks: a combinatorial multi-armed
bandit formulation,'' \emph{IEEE Symposium on International Dynamic
Spectrum Access Networks (DySPAN)}, April, 2010.

\bibitem{Liu:Zhao}
K.~Liu and Q.~Zhao, ``Decentralized multi-armed bandit with multiple
distributed players,'' \emph{Information Theory and Applications
Workshop (ITA)}, January, 2010.

\bibitem{Anandkumar:2010}
A.~Anandkumar, N.~Michael, and A.~K.~Tang. ``Opportunistic spectrum
access with multiple users: learning under competition,'' \emph{IEEE
International Conference on Computer Communications}, March, 2010.


\bibitem{Lai:Robbins}
T.~Lai and H.~Robbins, ``Asymptotically efficient adaptive
allocation rules,'' \emph{Advances in Applied Mathematics}, vol. 6,
no. 1, pp. 4-22, 1985.


\bibitem{Anantharam}
V.~Anantharam, P.~Varaiya, and J.~Walrand, ``Asymptotically
efficient allocation rules for the multiarmed bandit problem with
multiple plays-part I: IID rewards,'' \emph{IEEE Tran. on Auto.
Control}, vol. 32, no. 11, pp. 968-976, 1987.

\bibitem{Agrawal:1995}
R.~Agrawal, ``Sample mean based index policies with O(log n) regret
for the multi-armed bandit problem,'' \emph{Advances in Applied
Probability}, vol. 27, no. 4, pp. 1054-1078, 1995.

\bibitem{Auer:2002}
P.~Auer, N.~Cesa-Bianchi, and P.~Fischer, ``Finite-time analysis of
the multiarmed bandit problem,'' \emph{Machine Learning}, vol. 47,
no. 2, pp. 235-256, 2002.


\bibitem{Anantharam:1987}
V.~Anantharam, P.~Varaiya, and J .~Walrand, ``Asymptotically
efficient allocation rules for the multiarmed bandit problem with
multiple plays-part II: markovian rewards,'' \emph{IEEE Transactions
on Automatic Control}, vol. 32, no. 11, pp. 977-982, 1987.


\bibitem{Tekin:2010}
C.~Tekin, M.~Liu, ``Online algorithms for the multi-armed bandit
problem with markovian rewards,'' Allerton Conference, September,
2010.


\bibitem{Tekin:restless}
C. Tekin and M. Liu, ``Online learning in opportunistic spectrum
access: a restless bandit approach,'' \emph{IEEE International
Conference on Computer Communications (INFOCOM)}, April, 2011.

\bibitem{Qing:restless}
H. Liu, K. Liu, and Q. Zhao, ``Logarithmic weak regret of
non-Bayesian restless multi-armed bandit,'' \emph{IEEE International
Conference on Acoustics, Speech and Signal Processing (ICASSP)},
May, 2011.

\bibitem{Dai:restless}
W. Dai, Y. Gai, B. Krishnamachari and Q. Zhao, ``The non-Bayesian
restless multi-armed bandit: a case of near-logarithmic regret,''
\emph{IEEE International Conference on Acoustics, Speech and Signal
Processing (ICASSP)}, May, 2011.

\bibitem{Rusmevichientong}
P.~Rusmevichientong and J.~N.~Tsitsiklis, ``Linearly parameterized
bandits,'' \emph{Mathematics of Operations Research}, vol. 35 , no.
2, pp. 395-411, 2010


\bibitem{Kuhn:1955} H.~W.~Kuhn, ``The hungarian method for the assignment
problem,'' \emph{Naval Research Logistics Quarterly}, vol. 2, no. 1,
pp. 83-97, 1955.

\bibitem{Gillman:1998}
D.~Gillman, ``A chernoff bound for random walks on expander
graphs,'' \emph{SIAM Journal on Computing}, vol. 27, no. 4, pp.
1203-1220, 1998

\end{thebibliography}
\end{document}